\documentclass[a4paper,twocolumn]{article}

\usepackage[T1]{fontenc}
\usepackage[utf8]{inputenc}
\usepackage[english]{babel}
\usepackage{graphics}
\usepackage{epsfig}
\usepackage{amsmath} 
\usepackage{amssymb}
\usepackage{bbold}
\usepackage{url}
\usepackage{hyperref}
\usepackage{draftwatermark}
\usepackage{float}
\usepackage{pifont}
\usepackage[table]{xcolor}

\SetWatermarkText{Preprint}
\SetWatermarkScale{5}

\newcommand{\cmark}{\ding{51}}
\newcommand{\xmark}{\ding{55}}

\definecolor{lightgray}{gray}{0.8}

\DeclareMathOperator{\R}{\mathbb{R}}					
\DeclareMathOperator{\1}{\mathbb{1}}					
\DeclareMathOperator{\tr}{tr} 						
\DeclareMathOperator{\tov}{vec}
\DeclareMathOperator{\arsinh}{arsinh}

\newenvironment{definition}[1][Definition]{\begin{trivlist}
\item[\hskip \labelsep {\bfseries #1}]}{\end{trivlist}}

\title{\vspace{-3ex}A Unified Software Framework for Empirical Gramians}
\author{Christian Himpe\thanks{Contact: \href{mailto:christian.himpe@wwu.de}{\nolinkurl{christian.himpe@wwu.de}}, \href{mailto:mario.ohlberger@uni-muenster.de}{\nolinkurl{mario.ohlberger@uni-muenster.de}}, Institute for Computational and Applied Mathematics at the University of M\"unster, Einsteinstrasse~62, D-48149 M\"unster, Germany} \and Mario Ohlberger\footnotemark[1]}
\date{}

\begin{document}

\setlength{\parindent}{0pt}

\maketitle
\thispagestyle{empty}
\pagestyle{empty}


\begin{abstract}
\textbf{
A common approach in model reduction is balanced truncation, which is based on gramian matrices classifying certain attributes of states or parameters of a given dynamic system.
Initially restricted to linear systems, the empirical gramians not only extended this concept to nonlinear systems, but also provide a uniform computational method. 
This work introduces a unified software framework supplying routines for six types of empirical gramians.
The gramian types will be discussed and applied in a model reduction framework for multiple-input-multiple-output systems. 
}
\end{abstract}

~


\section{Introduction}
In a control system setting, balanced truncation is a well known technique for model reduction. 
Introduced by \cite{moore81}, gramian matrices were employed to determine controllability and observability of linear systems.
From these gramians a balancing transformation can be computed, enabling the truncation, for example, of states that are neither controllable nor observable.

With \cite{lall99}, empirical (controllability and observability) gramians were introduced, which correspond to the analytical gramians for linear systems, while extending the concept of system gramians to nonlinear systems which are generally given by:
\begin{align*}
  \dot{x}(t) &= f(x,u,p) \\
  y(t) &= g(x,u,p),
\end{align*}
with the system function $f$ and output function $g$ of states $x$, input $u$ and parameters $p$.
In the special case of an unparametrized linear system $f = Ax(t)+Bu(t)$ and $g = Cx(t)$,
these empirical gramians are computed by averaging simulations or experimental data with perturbations in inputs and initial states.

The \textbf{emgr} framework presented here encompasses six empirical gramians, namely the controllability, observability, cross, sensitivity, identifiability and joint gramian.
To adapt the computation of empirical gramians to the operating setting of the system, the initial state and the input are the main parameters which are perturbed by rotations and scaling.
The sets of rotations provided are $\{ \1 \}$ (unit matrix) and $\{ -\1,\1 \}$ (negative unit matrix and unit matrix).
Though these are very basic sets, and thus might not reflect all dynamics, especially with interrelated states and parameters, they allow a very efficient gramian assembly.
Scales may be freely chosen.
The subdivision of the scales may be linear, logarithmic or geometric.
Finally, there are several options to average against: the arithmetic average \cite{lall99}, the median, a steady state \cite{hahn02a}, and additionally, the principal components of the simulations or data via a proper orthogonal decomposition (POD). 

\section{Empirical Gramians}
Concerned with the reduction of states, the controllability, observability and cross gramian are presented next;
followed by the sensitivity, identifiability and joint gramian, which are used for parameter and combined reduction.
For the purpose of defining the gramians, a linear time-invariant control system is assumed:
\begin{align*}
  \dot{x}(t) &= Ax(t) + Bu(t) \\
  y(t) &= Cx(t),
\end{align*}
with the states $x \in \R^n$, control or input $u \in \R^m$, output $y \in \R^o$, system matrix $A \in \R^{n \times n}$, input matrix $B \in \R^{n \times m}$ and output matrix $C \in \R^{o \times n}$.

The necessary perturbations are given by six sets, of which $\{E_u,R_u,Q_u\}$ define the input perturbations, while sets $\{E_x,R_x,Q_x\}$ define the initial state perturbations:
\begin{align*}
 E_u &= \{ e_i \in \R^j ; \|e_i\| = 1 ; e_i e_{j \neq i} = 0; i=1,\ldots,m \} \\
 E_x &= \{ f_i \in \R^n ; \|f_i\| = 1 ; f_i f_{j \neq i} = 0; i=1,\ldots,n \} \\
 R_u &= \{ S_i \in \R^{j \times j} ; S_i^* S_i = \1 ; i = 1,\ldots,s \} \\
 R_x &= \{ T_i \in \R^{n \times n} ; T_i^* T_i = \1 ; i = 1,\ldots,t \} \\ 
 Q_u &= \{ c_i \in \R ; c_i > 0 ; i = 1,\ldots,q \} \\ 
 Q_x &= \{ d_i \in \R ; d_i > 0 ; i = 1,\ldots,r \}.
\end{align*}
These sets should correspond to the ranges in inputs and initial states the system is operating in.

\subsection{Controllability Gramian}
Controllability is a quantification of how well a state can be driven by input. 
Analytically, the controllability gramian is given by the smallest semi-positive definite solution of the Lyapunov equation: $AW_C + W_CA^T = -BB^T$.
If the underlying system is asymptotically stable, the controllability gramian can also be defined using the linear input-to-state map:
\begin{align*}
 W_C = \int_0^\infty e^{A \tau} B B^T e^{A^T \tau} d \tau.
\end{align*}
Following \cite{hahn02a}, the empirical controllability gramian is defined by:

\begin{definition} ~\\
 For sets $E_u$, $R_u$, $Q_u$, input $u(t)$ and input during the steady state $\bar{x}$, $\bar{u}$, the \textbf{empirical controllability gramian} is given by:
 \begin{align*}
  W_C &= \frac{1}{|Q_u| |R_u|} \sum_{h=1}^{|Q_u|} \sum_{i=1}^{|R_u|} \sum_{j=1}^{m} \frac{1}{c_h^2} \int_0^\infty \Psi^{hij}(t) d t \\
  \Psi^{hij}(t) &= (x^{hij}(t) - \bar{x})(x^{hij}(t) - \bar{x})^* \in \R^{n \times n}.
 \end{align*}
 With $x^{hij}$ being the states for the input configuration $u^{hij}(t) = c_h S_i e_j u(t) + \bar{u}$.
\end{definition}

Originally, in \cite{lall99}, $u(t)$ was restricted to $\delta(t)$, but extended in \cite{hahn02a} to arbitrary input under the name empirical covariance matrix.
$\bar{x}$ can be the arithmetic average, the median, the steady state or the principal components.
Restricting $R_u$ to $\{ -\1,\1 \}$  simplifies the input perturbation to:
\begin{align*}
 u^{hij}(t) =  -1^i q_h e_j u(t) + \bar{u}.
\end{align*}

\subsection{Observability Gramian}
Observability quantifies how well a change in a state is reflected by the output.
The analytical observability gramian is given by the smallest semi-positive definite solution of the Lyapunov equation: $AW_O + W_OA^T = -C^TC$.
Given an asymptotically stable underlying system the observability gramian can also be defined using the state-to-output map:
\begin{align*}
 W_O = \int_0^\infty e^{A^T \tau} C^T C e^{A \tau} d \tau.
\end{align*}
The empirical observability gramian is defined as described in \cite{lall99} and \cite{hahn02a}:

\begin{definition} ~\\
For sets $E_x$, $R_x$, $Q_x$ and output $y$ during the steady state $\bar{x}$, $\bar{y}$, the \textbf{empirical observability gramian} is given by:
 \begin{align} \label{obs}
  W_O &= \frac{1}{|Q_x| |R_x|} \sum_{k=1}^{|Q_x|} \sum_{l=1}^{|R_x|} \frac{1}{d_k^2} T_l \int_0^\infty \Psi^{kl}(t) d t \; T_l^* \quad \\ \notag
  \Psi_{ab}^{kl} &= (y^{kla}(t) - \bar{y})^* (y^{klb}(t) - \bar{y}) \in \R.
 \end{align}
With $y^{kla}$ being the systems output for the initial state configuration $x^{kla}_0 = d_k S_l f_a + \bar{x}$.
\end{definition}

$\bar{y}$ can be the arithmetic average, the median, the steady state output or the principal components.
Restricting $R_x$ to $\{ -\1,\1 \}$ simplifies equation~(\ref{obs}) to:
\begin{align*}
 W_O &= \frac{1}{|Q_x| |R_x|} \sum_{k=1}^{|Q_x|} \sum_{l=1}^{|R_x|} \frac{1}{d_k^2} \int_0^\infty \Psi^{kl}(t) d t,
\end{align*}
and the initial state perturbation to:
\begin{align*}
 x^{kla}_0 = -1^l d_k f_a + \bar{x}.
\end{align*}

\subsection{Cross Gramian}
The cross gramian \cite{antoulas05} makes a combined statement about the controllability and observability, given the system has the same number of inputs and outputs.
If the system is also symmetric, meaning the system transfer function is symmetric, then the absolute value of this gramians' eigenvalues equal the Hankel singular values.
It is originally computed as the smallest semi-positive definite solution of the Sylvester equation: $AW_X + W_XA^T = -BC$.
The cross gramian can also be defined using the input-to-state and state-to-output maps, if the underlying system is asymptotically stable:
\begin{align*}
 W_X = \int_0^\infty e^{A \tau} B C e^{A \tau} d \tau
\end{align*}
The empirical cross gramian has been introduced in \cite{streif06} for SISO systems and was extended to MIMO systems in \cite{himpe13a}. 

\begin{definition} ~\\
For sets $E_u$, $E_x$, $R_u$, $R_x$, $Q_u$, $Q_x$, input $\bar{u}$ during steady state $\bar{x}$ with output $\bar{y}$, the \textbf{empirical cross gramian} is given by:
 \begin{align} \label{crs}
  W_X &= \frac{1}{|Q_u| |R_u| m |Q_x| |R_x|} \sum_{h=1}^{|Q_u|} \sum_{i=1}^{|R_u|} \sum_{j=1}^m \sum_{k=1}^{|Q_x|} \sum_{l=1}^{|R_x|} \notag \\ & \frac{1}{c_h d_k} \int_0^\infty T_l \Psi^{hijkl}(t) T_l^* d t \notag \\
  &\Psi_{ab}^{hijkl}(t)  = f_b^* T_k^* \Delta x^{hij}(t) e_i^* S_h^* \Delta y^{kla}(t) \\
  &\Delta x^{hij}(t) = (x^{hij}(t)-\bar{x} ) \notag \\
  &\Delta y^{kla}(t) = (y^{kla}(t)-\bar{y} ). \notag 
 \end{align}
With $x^{hij}$ and $y^{kla}$ being the states and output for the input $u^{hij}(t) = c_h S_i e_j u(t) + \bar{u}$ and initial state $x^{kla}_0 = d_k T_l f_a + \bar{x}$ respectively.
\end{definition}

$\bar{x}$ and $\bar{y}$ can be the arithmetic average, the median, the steady state or the principal components of the output.
Again, restricting $R_u$ and $R_x$ to $\{ -\1,\1 \}$ simplifies equation~(\ref{crs}) to:
\begin{align*}
 &\Psi_{ab}^{hijkl}  = (-1)^{i+l} f_b^* \Delta x^{hij}(t) e_j^* \Delta y^{kla}(t),
\end{align*}
as well as input and initial state perturbation to:
\begin{align*}
 u^{hij}(t) &= -1^i c_h e_j u(t) + \bar{u} \\
 x^{kla}_0 &= -1^l d_k f_a + \bar{x}.
\end{align*}

\subsection{Sensitivity Gramian}
The sensitivity gramian allows controllability based parameter reduction and identification.
It is based on \cite{sun06} and aimed for models that can be partitioned as follows:
\begin{align*}
 \dot{x} &= f(x,u,p) = f(x,u) + \sum_{k=1}^P f(x,p_k).
\end{align*}
The parameters $p \in \R^P$ are handled here as additional inputs.
All summands of the partitioned system are treated as independent subsystems and thus a controllability gramian for each subsystem can be computed.
Each parameters controllability is encoded in the sum of singular values of the associated sub-controllability gramian $W_{C,k}$.
The sensitivity gramian is now given by the diagonal matrix, with each diagonal element being the trace of a sub-controllability gramian:
\begin{align*}
 W_S = \sum_{k=1}^P \tr(W_{C,k}) (\delta_{ij})_{i=j=k}.
\end{align*}
The controllability and thus identifiability of each parameter is then given by the corresponding diagonal entry of the sensitivity gramian $W_S$.
For partionable linear systems, the sum of all subsystems controllability gramians $W_{C,k}$ and the parameter free subsystems gramian $W_{C,0}$ equals the usual controllability gramian \cite{sun06}: 
\begin{align*}
 W_C = W_{C,0} + \sum_{k=1}^P W_{C,k}.
\end{align*}
The sensitivity gramian can be applied to non-partitionable models with reduced accuracy.
 
\subsection{Identifiability Gramian}
The identifiability gramian enables observability based parameter identification and consequently parameter reduction.
As described in \cite{geffen08}, the dynamic systems states are augmented with as many states as parameters that are constant over time, and have the initial value of the (prior) parameter value.
\begin{align*}
 \dot{\breve{x}} = \begin{pmatrix} \dot{x} \\ \dot{p} \end{pmatrix} &= \begin{pmatrix} f(x,u,p) \\ 0 \end{pmatrix} \\ 
 y &= g(x,u,p)
\end{align*}
The observability gramian of this augmented system holds the observability information of states and parameters.
To extract the parameter specific observability, the Schur-complement can be applied to the augmented observability gramian:
\begin{align*}
 W &= \begin{pmatrix} W_O & \vline & W_Q \\ \hline W_Q^* & \vline & W_P \end{pmatrix}\\
 \Rightarrow W_I &= W_P - W_Q^* {W_O}^{-1} W_Q.
\end{align*}

\subsection{Joint Gramian}
Based on the identifiability gramian procedure, the cross gramian can be employed for a concurrent state and parameter reduction (see \cite{himpe13a}).
Not only augmenting the states with constant parameter states as before, but also adding as many inputs $v$ and outputs as parameters, acting via identity on the augmented states and augmented outputs
\begin{align*}
 \dot{\breve{x}} &= \begin{pmatrix} \dot{x} \\ \dot{p} \end{pmatrix} = \begin{pmatrix} f(x,u,p) \\ v \end{pmatrix} \\ 
 y &= \begin{pmatrix} g(x,u,p) \\ \1 \end{pmatrix},
\end{align*}
to preserve symmetry.
The cross gramian of this special augmented system, similar to the identifiability gramian, holds the cross gramian of the original system as well as a cross identifiability gramian $W_{\ddot{I}}$ which can be extracted with the Schur-complement from the joint gramian:
\begin{align*}
 W_J &= \begin{pmatrix} W_X & \vline & W_Q \\ \hline W_Q^* & \vline & W_P \end{pmatrix}\\
 \Rightarrow W_{\ddot{I}} &= W_P - W_Q^* {W_X}^{-1} W_Q.
\end{align*}

\section{Implementation}
The \textbf{emgr} software framework presented here provides a uniform interface to compute all six empirical gramians and is given by:
\begin{align*}
 \mathtt{W = emgr(f,g,q,p,t,w,v_{cfg},u,u_s,x_s,u_m,x_m,y_d);}
\end{align*}
with $\mathtt{f}$ and $\mathtt{g}$ being handles to the system and the output function; both requiring the signature $\mathtt{f(x,u,p)}$ and $\mathtt{g(x,u,p)}$.
$\mathtt{q}$ is a vector defining the systems number of inputs, states and outputs.
$\mathtt{p}$ holds the parameters.
$\mathtt{t}$ is a  three component vector containing start time, time step and stop time.
$\mathtt{w}$ is a character setting the gramian type; for an overview on the applicability of gramian types see table~\ref{tab}.

Following arguments are optional.
The ten component vector $\mathtt{v_{cfg}}$ configures the available options, including averaging types, input and state scale subdivisions and perturbation rotations.
$\mathtt{u}$ provides the input to $f$ and $g$, while $\mathtt{u_s}$ and $\mathtt{x_s}$ set steady input and steady state.
$\mathtt{u_m}$ and $\mathtt{x_m}$ define the scales of the perturbation.
Lastly, $\mathtt{y_d}$ allows to pass experimental data to be used instead of generated snapshots.

The parameter reducing empirical gramians (sensitivity, identifiability and joint) are an encapsulation of the state reducing empirical gramians (controllability, observability and cross).
Computation of the latter is extensively vectorized, exploiting the gramian matrix assembly format.
In example, the empirical observability gramian assembly from equation~(\ref{obs}) can be computationally simplified to:
\begin{align*}
 \Psi_{ab}^{kl} = (y^{kla}(t) - \bar{y})^* (y^{klb}(t) - \bar{y}) \\
 \to \begin{cases} \psi_{a}^{kl}(t) = (y^{kla}(t) - \bar{y}) \\ \Psi_{ab}^{kl} = \tov(\psi_{a}^{kl})^*\tov(\psi_{b}^{kl}). \end{cases}
\end{align*}
Computation of empirical gramians using \textbf{emgr} is very portable, since only basic vector and matrix operations are required.
Necessary integrations, meaning simulations for given inputs or initial states, are accomplished either by the first order Euler's method, second order Adams-Bashfort method or second order leapfrog method.
The empirical gramian framework \textbf{emgr} as well as the following experiments are released under an open source license, are compatible with OCTAVE and MATLAB, and can be found at \mbox{\url{http://gramian.de}} or at the \mbox{\href{http://orms.mfo.de/project?id=345}{Oberwolfach References on Mathematical Software}}.

\begin{table}[t]\footnotesize
\rowcolors{3}{white}{lightgray}
\begin{tabular}{c|c|c|c}
 \textbf{Gramian} & \textbf{State}    & \textbf{Parameter} & \textbf{Combined} \\
  \textbf{Type}   & \textbf{Reduction} & \textbf{Reduction} & \textbf{Reduction} \\
\hline
 $W_C$ & (\cmark) & \xmark & \xmark \\
 $W_O$ & (\cmark) & \xmark & \xmark \\
 $W_C \& W_O$ & \cmark & \xmark & \xmark \\
 $W_X$ &  \cmark  & (\cmark) & \xmark \\
 $W_S$ & (\cmark) &  \cmark & \xmark \\ 
 $W_I$ & (\cmark) &  \cmark & (\cmark) \\
 $W_J$ & (\cmark) &  \cmark & \cmark \\
\end{tabular}
\caption{Empirical gramian application matrix.}
\label{tab}
\end{table}

\section{Numerical Experiments}
To demonstrate the various empirical gramians, computed by the \textbf{emgr} framework, a symmetric nonlinear MIMO system with one hundred states, ten inputs and ten outputs is employed. 
The system matrix is generated randomly with ensured stability and symmetry; the input matrix $B$ is also a random matrix and the output matrix is given by $C = B^T$.
Furthermore a random, but element-wise, parametrized source term $E_p$ of dimension $n=100$ parameters is added.
Input is applied through a delta impulse.
\begin{align*}
 f(x,u,p) &= \dot{x} = A \arsinh(x) + Bu + E_p\\
 g(x) &= y = Cx 
\end{align*}

\begin{figure}[h!]
 \centerline{\includegraphics[scale=0.45]{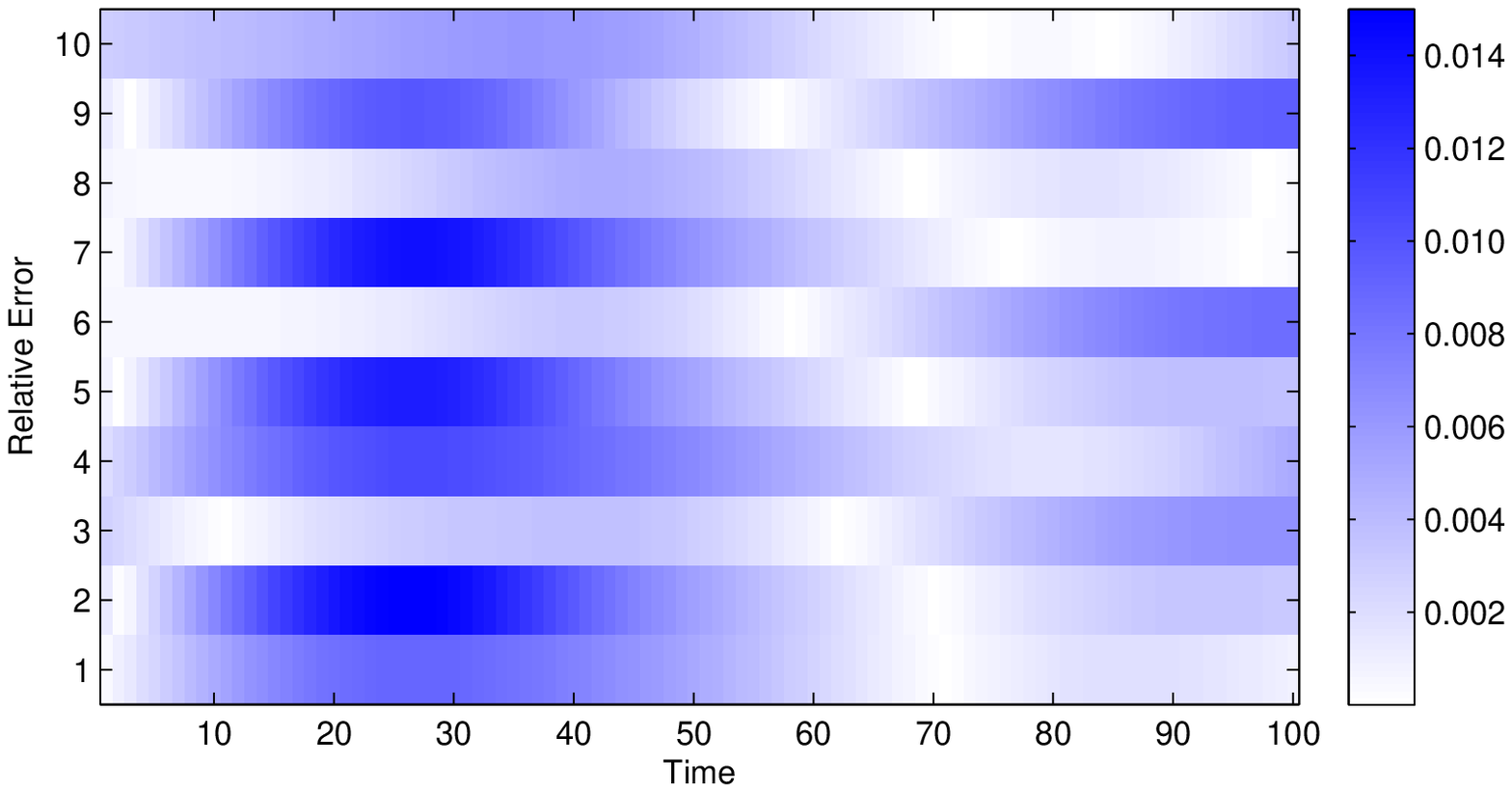}}
 \caption{\textit{Relative error in reduced system output by balanced truncation using the empirical controllability gramian and the empirical observability gramian} \small{\texttt{WC~=~emgr(f,g,[10,100,10],p,[0,0.01,1],'c');} \texttt{WO~=~emgr(f,g,[10,100,10],p,[0,0.01,1],'o');}} } 
 \label{wcwo}
\end{figure}

First, a state reduction, using the empirical controllability gramian and the empirical observability gramian, through balanced truncation is performed in figure~\ref{wcwo}, reducing the number of states to the number of outputs.
Balanced truncation as a classic approach in model order reduction will be used as a baseline, to which the following methods will be compared.

\begin{figure}[h!]
 \centerline{\includegraphics[scale=0.45]{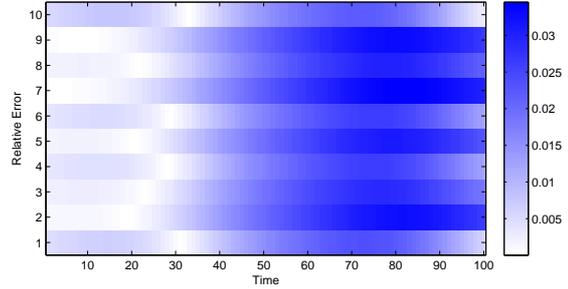}}
 \centering
 \caption{\textit{Relative error in reduced system output by truncation using the empirical cross gramian} \small{\texttt{WX~=~emgr(f,g,[10,100,10],p,[0,0.01,1],'x');}}}
 \label{wx}
\end{figure}

Next, a state reduction by direct truncation employing the empirical cross gramian is demonstrated in figure~\ref{wx}; again reducing the number of states to ten.
The state reduction via direct truncation of the cross gramian has about the same error, but requires only half of the reduction time, since only one gramian and no balancing is required. 

\begin{figure}[h!]
 \centerline{\includegraphics[scale=0.45]{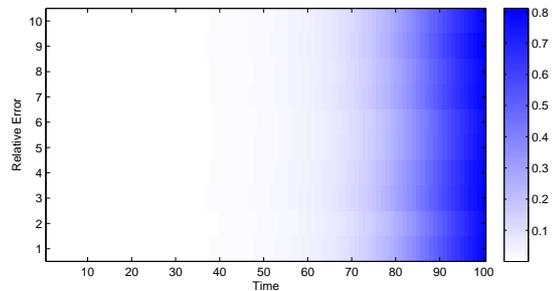}}
 \centering
 \caption{\textit{Relative error in system with reduced source term; reduction by truncation using the empirical sensitivity gramian} \small{\texttt{WS~=~emgr(f,g,[10,100,10],p,[0,0.01,1],'s');}} }
 \label{ws}
\end{figure}

The empirical sensitivity gramian can be applied if the underlying system can be partitioned such that $f(x,u,p) = f(x,u) + f(x,p)$. 
To be able to use it in this setting, the parametrized source term is reduced to the number of outputs in figure~\ref{ws}.
The sensitivity gramian is the fastest parameter reduction method, but has a high relative error in outputs.

\begin{figure}[h!]
 \centerline{\includegraphics[scale=0.45]{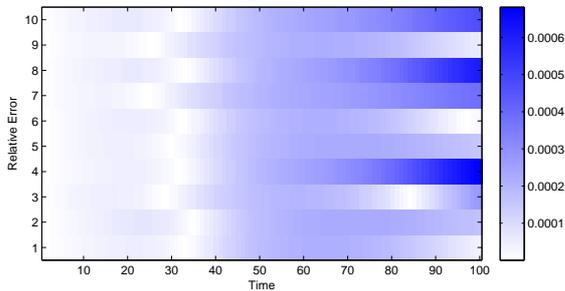}}
 \centering
 \caption{\textit{Relative error in reduced with reduced source term; reduction by truncation using the empirical identifiability gramian} \small{\texttt{WI~=~emgr(f,g,[10,100,10],p,[0,0.01,1],'i');}} }
 \label{wi}
\end{figure}

Since the parameters of the source term are reduced, the cumulative effects in the original system are the origin of the increasing error over time.
Next, the parametrized source term is reduced by the empirical identifiability gramian in figure~\ref{wi}.
Taking five times as long for the parameter reduction, the identifiability gramian is about two orders of magnitude more accurate.

\begin{figure}[h!]
 \centerline{\includegraphics[scale=0.45]{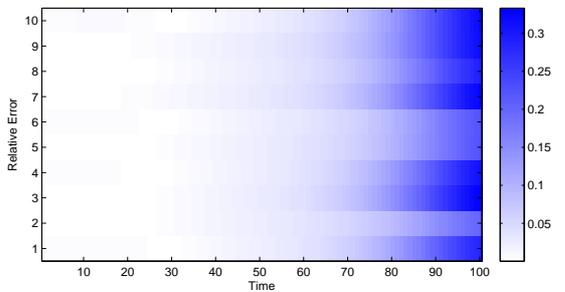}}
 \centering
 \caption{\textit{Relative error in system output with combined state and parameter reduction; reduction by truncation of parameters and states using the empirical joint gramian} \small{\texttt{WJ~=~emgr(f,g,[10,100,10],p,[0,0.01,1],'j');}} }
 \label{wj}
\end{figure}

Finally, in figure~\ref{wj}, the same system undergoes a combined state and parameter reduction using the empirical joint gramian.
Though with longest total duration, the joint gramian is the only gramian allowing direct, balancing-free combined reduction of state and parameter space with a comparable relative error. 
This combined reduction generates a reduced-order model, of which the relative error is comparable to the other reduced models. 

\section{Future Work}
The \textbf{emgr} framework already allows a wide range of computations of empirical gramians for state or parameter reduction.
Apart from model order reduction, the empirical gramians can be employed for system identification tasks, like parameter identification or sensitivity analysis as well as decentralized control, nonlinearity measurement and uncertainty quantification. 

Further work will enhance the flexibility, while keeping the interface as simple as possible.
Following \cite{geffen08}, allowing factorial designs will greatly enlarge the field of application.
Finally, extending the use of the cross gramian (and thus the joint gramian) to non-symmetric systems \cite{antoulas05} will enable a combined state and parameter reduction for general linear and nonlinear models without balancing.


\section*{Acknowledgement}

We acknowledge support by the Deutsche Forschungsgemeinschaft and Open Access Publication Fund of the University of M\"unster

\end{document}